\documentclass[12pt]{article}

\newtheorem{definition}{Definition}
\newtheorem{theorem}{Theorem}
\newtheorem{proposition}{Proposition}
\newtheorem{lemma}{Lemma}

\newtheorem{corollary}{Corollary}

\usepackage{amssymb}
\usepackage{latexsym}

\begin{document}

\begin{flushleft}
{\bf Generalized Shannon inequalities based on Tsallis relative operator entropy}
\end{flushleft}
\vspace{0.2cm}
\begin{flushleft}
Kenjiro Yanagi\footnote{E-mail:yanagi@yamaguchi-u.ac.jp} \\
Department of Applied Science \\
Faculty of Engineering \\
Yamaguchi University \\
Tokiwadai 2-16-1 \\
Ube 755-0811, Japan, 
\end{flushleft}
\begin{flushleft}
Ken Kuriyama \footnote{E-mail:kuriyama@yamaguchi-u.ac.jp} \\
Department of Applied Science \\
Faculty of Engineering \\
Yamaguchi University \\
Tokiwadai 2-16-1 \\
Ube 755-0811, Japan 
\end{flushleft}
\begin{flushleft}
and 
\end{flushleft}
\begin{flushleft}
Shigeru Furuichi \footnote{E-mail:furuichi@ed.yama.tus.ac.jp} \\
Department of Electronics and Computer Science \\
Faculty of Science and Engineering \\
Tokyo University of Science, Yamaguchi \\
Daigaku-Dori 1-1-1 \\
Onoda 756-0884, Japan
\end{flushleft}
\vspace{0.5cm}
\begin{flushleft}
{\sc Abstract}  
\end{flushleft}

Tsallis relative operator entropy is defined and then its properties are
given.  Shannon inequality and its reverse one in Hilbert space operators derived by T.Furuta \cite{Fu:par}
 are extended in terms of the parameter of the Tsallis relative operator entropy. Moreover the generalized Tsallis relative operator entropy is introduced and then
several operator inequalities are derived.

\vspace{3mm}

\vspace{3mm}


%
\vspace{0.5cm}
\begin{flushleft}
{\sc Introduction}
\end{flushleft}

Tsallis entropy 
$$
S_q(X) = -\sum_{x} p(x)^q \ln_q p(x)
$$
was defined in \cite{Ts:rel} for the probability distribution $p(x)$, 
where $q$-logarithm function is defined by 
$\ln_q(x)\equiv \frac{x^{1-q}-1}{1-q}$ 
for any nonnegative real numbers $x$ and $q \neq 1$. 
It is easily seen that Tsallis entropy is one parameter extension of 
Shannon entropy $S_1(X) \equiv -\sum_x p(x) \log p(x)$ and converges to it as 
$q \to 1$.  
The study based on Tsallis type entropies has been developed in mainly 
statistical physics \cite{AO:re1}. 
In the recent work \cite{Ab:re1}, Tsallis type relative entropy in quantum 
system,  defined by 
\begin{equation}\label{dq}
D_q(\rho\vert \sigma) \equiv \frac{1}{1-q} \left[1- Tr(\rho^q\sigma^{1-q})\right]
\end{equation}
for two density operators $\rho$ and $\sigma$ (i.e., positive operators with unit trace) 
and $0\leq q<1$,  was investigated. 
 
  On the other hand, the relative operator entropy was defined by 
  J.I.Fujii and E.Kamei in \cite{FuKa:rel}. 
Many important results in operator theory and information theory have been 
published in the relation to Golden-Thonmpson inequality \cite{HP,BPL}. 
We are interested in not only the properties of the Tsallis type relative 
entropy  but also the properties before taking a trace, namely, 
Tsallis type relative operator entropy which is a parametric extension of 
the relative operator entropy. 
In this paper, we define the Tsallis relative operator entropy and then show 
some properties of Tsallis relative operator entropy. To this end, we slightly change the 
parameter $q$ in Eq.(\ref{dq}) to $\lambda$ in our definition which will be 
appeared in the following section. Moreover, in order to make our definition 
correspond to the definition of the relative operator entropy defined in 
\cite{FuKa:rel}, we change the sign of the original Tsallis relative entropy. 

%
\section{Tsallis relative entropy}

As mentioned above, we adopt the slightly modified definition of the Tsallis relative entropy in the following. 

\begin{definition}
Let $a = \{a_1,a_2,\ldots,a_n \}$ and $b = \{b_1,b_2,\ldots,b_n \}$ 
be two probability vectors satisfying $a_j, b_j > 0$. Then for 
$0 < \lambda \leq 1$
\begin{equation}\label{our}
S_{\lambda}(a|b) = \frac{\sum_{j=1}^n a_j^{1-\lambda}b_j^{\lambda}-1}{\lambda}
\end{equation}
is called Tsallis relative entropy between $a$ and $b$. 
\label{def:definition1}
\end{definition}

We should note that the Tsallis relative entropy is usually defined by
\begin{equation} \label{usually}
 D_q(a\vert b) = \frac{1-\sum_{j=1}^n a_j^q b_j^{1-q}}{1-q},
\end{equation}
with a parameter $q \geq 0$ in the field of statisitical physics \cite{AO:re1}. There is the relation between them 
such that $S_{\lambda}(a \vert b)=-D_{1-q}(a \vert b)$.
However, in this paper, we adopt the definition of Eq.(\ref{our}) in stead of Eq.(\ref{usually}),
 to study of the properties of the parametrically extended relative operator entropy as a series of the study of 
the relative operator entropy from the operator theoretical point of view. 
The opposite sign between the relative entropy defined by Umegaki \cite{Um} and the relative operator entropy led
 us to define the Tsallis relative operator entropy in the above.

Tsallis relative entropy defined in Eq.(\ref{our}) has the following 
properties. 

\begin{proposition}
We have the following {\rm (1)} and {\rm (2)}. 
\begin{description}
\item[(1)] $\displaystyle{S_{\lambda}(a|b) \geq \sum_{j=1}^n a_j \log \frac{b_j}{a_j}}$ for $0 < \lambda \leq 1$.
\item[(2)] $\displaystyle{\lim_{\lambda \to 0} S_{\lambda}(a|b) = \sum_{j=1}^n a_j \log \frac{b_j}{a_j}}$.
\end{description}
\label{prop:proposition1}
\end{proposition}

\begin{flushleft}
{\bf Proof.}  (1)  Since $t^{\lambda}-1 \geq \log t^{\lambda}$, we have 
\end{flushleft}
$$
\frac{\sum_{j=1}^n a_j^{1-\lambda}b_j^{\lambda}-1}{\lambda} = \sum_{j=1}^n a_j\frac{(\frac{b_j}{a_j})^{\lambda}-1}{\lambda} \geq \sum_{j=1}^n a_j \log \frac{b_j}{a_j}.
$$

\noindent
(2) 
\begin{eqnarray*}
\lim_{\lambda \to 0} \frac{\sum_{j=1}^n a_j^{1-\lambda}b_j^{\lambda}-1}{\lambda} & = & \lim_{\lambda \to 0} \frac{\sum_{j=1}^n a_j(\frac{b_j}{a_j})^{\lambda}-1}{\lambda} \\
& = & \sum_{j=1}^n a_j(\frac{b_j}{a_j})^{\lambda}\log \frac{b_j}{a_j} |_{\lambda = 0} \\
& = & \sum_{j=1}^n a_j \log \frac{b_j}{a_j}. 
\end{eqnarray*}
\ \hfill q.e.d.

\begin{proposition}
$$
0 \geq S_{\lambda}(a|b) \geq \frac{(\sum_{j=1}^n \frac{a_j^2}{b_j})^{-\lambda}-1}{\lambda}. 
$$
\label{prop:proposition2}
\end{proposition}

\begin{flushleft}
{\bf Proof.} Since 
\end{flushleft}
$$
\sum_{j=1}^n a_j^{1-\lambda}b_j^{\lambda} \leq \sum_{j=1}^n \{(1-\lambda)a_j+\lambda b_j \} = (1-\lambda)\sum_{j=1}^n a_j + \lambda \sum_{j=1}^n b_j = 1, 
$$
we have 
$$
S_{\lambda}(a|b) \leq 0.
$$
We also give another inequality. Since 
$$
\sum_{j=1}^n a_j^{1-\lambda}b_j^{\lambda} = \sum_{j=1}^n (\frac{a_j}{b_j})^{-\lambda}a_j \geq \prod_{j=1}^n (\frac{a_j}{b_j})^{-\lambda a_j} = (\prod_{j=1}^n (\frac{a_j}{b_j})^{a_j})^{-\lambda} \geq (\sum_{j=1}^n \frac{a_j^2}{b_j})^{-\lambda}, 
$$
we have 
$$
S_{\lambda}(a|b) \geq \frac{(\sum_{j=1}^n \frac{a_j^2}{b_j})^{-\lambda}-1}{\lambda}. 
$$
\ \hfill q.e.d. 

%
\section{Tsallis relative operator entropy}

A bounded linear operator $T$ on a Hilbert space $H$ is said to be positive 
(denoted by $T \geq 0$) if $(Tx,x) \geq 0$ for all $x \in H$ and also an 
operator $T$ is said to be strictly positive (denoted by $T > 0$) if $T$ is 
invertible and positive.  We define Tsallis relative operator entropy 
in the following. 

\begin{definition}
For $A > 0, B > 0$ and $0 < \lambda \leq 1$, 
$$
T_{\lambda}(A|B) = \frac{A^{1/2}(A^{-1/2}BA^{-1/2})^{\lambda}A^{1/2}-A}{\lambda}$$
is called Tsallis relative operator entropy between $A$ and $B$. 
\label{def:definition2}
\end{definition}

In this section we give the Shannon type operator inequality and its reverse one satisfied 
by Tsallis relative operator entropy. 

\begin{theorem}
Let $\{ A_1,A_2,\ldots,A_n \}$ and $\{B_1,B_2,\ldots,B_n \}$ be two sequences 
of strictly positive operators on a Hilbert space $H$. If 
$\displaystyle{\sum_{j=1}^n A_j = \sum_{j=1}^n B_j = I}$, then 
$$
0 \geq \sum_{j=1}^n T_{\lambda}(A_j|B_j) \geq \frac{(\sum_{j=1}^n A_jB_j^{-1}A_j)^{-\lambda}-I}{\lambda}. 
$$
\label{th:theorem1}
\end{theorem}

We need a lemma before we prove the main theorem. 

\begin{lemma}
For fixed $t > 0$, an inequality of 
$\lambda \; (0 < \lambda \leq 1)$ holds. 
$$
\frac{t^{\lambda}-1}{\lambda} \leq t-1.
$$
\label{lem:lemma1}
\end{lemma}

\begin{flushleft}
{\bf Proof.}  If $t = 1$, then it is clear. If $t \neq 1$, then 
we put $F(\lambda) = \lambda(t-1)-t^{\lambda}+1$. Then we have 
$F^{'}(\lambda) = t-1-t^{\lambda}\log t$ and 
$F^{''}(\lambda) = -t^{\lambda}(\log t)^2 < 0$. Hence $F(\lambda)$ is concave. 
Since $F(0) = F(1) = 0$, we have the result. 
\ \hfill q.e.d.
\end{flushleft}

\begin{flushleft}
{\bf Proof of Theorem \ref{th:theorem1}.} It follows from Lemma \ref{lem:lemma1} that 
\end{flushleft}
\begin{eqnarray*}
\frac{A^{1/2}(A^{-1/2}BA^{-1/2})^{\lambda}A^{1/2}-A}{\lambda} & = & 
A^{1/2} \frac{(A^{-1/2}BA^{-1/2})^{\lambda}-I}{\lambda} A^{1/2} \\
& \leq & A^{1/2}(A^{-1/2}BA^{-1/2}-I)A^{1/2} \\
& = & B-A, 
\end{eqnarray*}
where $A > 0, B > 0$ and $0 < \lambda \leq 1$. 
Then we have 
\begin{eqnarray*}
\sum_{j=1}^n T_{\lambda}(A_j|B_j) & = & \sum_{j=1}^n \frac{A_j^{1/2}(A_j^{-1/2}B_jA_j^{-1/2})^{\lambda}A_j^{1/2}-A_j}{\lambda} \\
& \leq & \sum_{j=1}^n (B_j - A_j) = 0.
\end{eqnarray*}
We also prove another inequality. We apply Proposition 3.1 of 
Furuta \cite{Fu:par} by putting $f(x) = -x^{-\lambda}$, $C_j = A_j^{1/2}$ and 
$X_j = A_j^{1/2}B_j^{-1}A_j^{1/2}$. Then 
$$
-(\sum_{j=1}^n A_j^{1/2}(A^{1/2}B_j^{-1}A_j^{1/2})A_j^{1/2})^{-\lambda} \geq -\sum_{j=1}^n A_j^{1/2}(A_j^{1/2}B_j^{-1}A_j^{1/2})^{-\lambda}A_j^{1/2}.
$$
Hence 
$$
(\sum_{j=1}^n A_jB_j^{-1}A_j)^{-\lambda} \leq \sum_{j=1}^n A_j^{1/2}(A_j^{-1/2}B_jA_j^{-1/2})^{\lambda}A_j^{1/2}.
$$
Then we complete the proof. \ \hfill q.e.d. \\

We also obtain the operator version of the Shannon inequality and reverse one given by Furuta \cite{Fu:par} 
as a corollary of Theorem \ref{th:theorem1} in the following.

\begin{corollary}[Furuta \cite{Fu:par}]
Let $\{ A_1,A_2,\ldots,A_n \}$ and $\{B_1,B_2,\ldots,B_n \}$ be two sequences 
of strictly positive operators on a Hilbert space $H$. If 
$\displaystyle{\sum_{j=1}^n A_j = \sum_{j=1}^n B_j = I}$, then 
$$
0 \geq \sum_{j=1}^n A_j^{1/2}(\log A_j^{-1/2}B_jA_j^{-1/2})A_j^{1/2} \geq -\log [\sum_{j=1}^n A_jB_j^{-1}A_j]. 
$$
\label{cor:corollary1}
\end{corollary}

We need the following lemma to prove it.

\begin{lemma}
For $0 < \lambda < 1,  0 < \alpha < \beta$, we have the following {\rm (1)} 
and {\rm (2)}. 
\begin{description}
\item[(1)]  $\displaystyle{\lim_{\lambda \to +0}\frac{t^{\lambda}-1}{\lambda} = \log t}$  uniformly on $[\alpha,\beta]$. 
\item[(2)]  $\displaystyle{\lim_{\lambda \to +0}\frac{t^{-\lambda}-1}{\lambda} = -\log t}$ uniformly on $[\alpha,\beta]$.
\end{description}
\label{lem:lemma2}
\end{lemma}

\begin{flushleft}
{\bf Proof.} We prove it by using Dini's theorem. \ \hfill q.e.d.
\end{flushleft}

\begin{flushleft}
{\bf Proof of Corollary \ref{cor:corollary1}.}  By (1) of Lemma \ref{lem:lemma2}, 
\end{flushleft}
$$
\lim_{\lambda \to +0} \frac{(A^{-1/2}BA^{-1/2})^{\lambda}-I}{\lambda} = \log A^{-1/2}BA^{-1/2}, 
$$
where the limit is taken in operator norm. Then 
\begin{eqnarray*}
\lim_{\lambda \to +0} \sum_{j=1}^n T_{\lambda}(A_j|B_j) & = & \lim_{\lambda \to +0} \sum_{j=1}^n \frac{A_j^{1/2}(A_j^{-1/2}B_jA_j^{-1/2})^{\lambda}A_j^{1/2}-A_j}{\lambda} \\
& = & \sum_{j=1}^n A_j^{1/2} (\log A_j^{-1/2}B_jA_j^{-1/2})A_j^{1/2}. 
\end{eqnarray*}
On the other hand, by (2) of Lemma \ref{lem:lemma2},  we have 
$$
\lim_{\lambda \to +0}\frac{(\sum_{j=1}^nA_jB_j^{-1}A_j)^{-\lambda}-I}{\lambda}=-\log[\sum_{j=1}^nA_jB_j^{-1}A_j]$$
Therefore Theorem \ref{th:theorem1} ensures 
$$
0 \geq \sum_{j=1}^nA_j^{1/2}(\log A_j^{-1/2}B_jA_j^{-1/2})A_j^{1/2}\geq -\log[\sum_{j=1}^nA_jB_j^{-1}A_j].$$

\ \hfill q.e.d.

Actually the above Corollary \ref{cor:corollary1} is a part of the Corollary 2.4 in \cite{Fu:par}. 
We will generalize our Tsallis relative operator entropy and derive some generalized operator inequalities 
by the different way from \cite{Fu:par} in the following section.

%
\section{Generalized Tsallis relative operator entropy}

We remind of the relative operator entropy and its related operator entropy. 

\begin{definition}
For $A > 0, B > 0$ 
$$
S(A|B) = A^{1/2}(\log A^{-1/2}BA^{-1/2})A^{1/2}
$$
is called relative operator entropy between $A$ and $B$. It was defined by 
Fujii and Kamei {\rm \cite{FuKa:rel}} originally. For $A > 0, B > 0$ and 
$\lambda \in \mathbb{R}$, the generalized relative operator entropy was defined by Furuta in \cite{Fu:par} 
$$
S_{\lambda}(A|B) = A^{1/2}(A^{-1/2}BA^{-1/2})^{\lambda}(\log A^{-1/2}BA^{-1/2})A^{1/2}
$$
and 
$$
A \natural_{\lambda} B = A^{1/2}(A^{-1/2}BA^{-1/2})^{\lambda}A^{1/2}.
$$
In particular we remark that $S_0(A|B) = S(A|B), A \natural_0 B = A$ and $A \natural_1 B = B$.
\label{def:definition3}
\end{definition}

We generalize the definition of the Tsallis relative operator entropy. 

\begin{definition}
For $A > 0, B > 0$, $\lambda, \mu \in \mathbb{R}, \lambda \neq 0$ and 
$k \in \mathbb{Z}$, 
$$
\tilde{T}_{\mu,k,\lambda}(A|B) = \frac{A \natural_{\mu+k\lambda} B-A \natural_{\mu+(k-1)\lambda} B}{\lambda}
$$
is called generalized Tsallis relative operator entropy. In particular we remark that for $\lambda \neq 0$ 
$$
\tilde{T}_{0,1,\lambda}(A|B) = \frac{A \natural_{\lambda} B-A \natural_{0} B}{\lambda} = \frac{A^{1/2}(A^{-1/2}BA^{-1/2})^{\lambda}A^{1/2}-A}{\lambda} = T_{\lambda}(A|B).
$$
\label{def:definition4}
\end{definition}

We state the relationship among $S_{\mu \pm k \lambda}(A|B), S_{\mu \pm (k+1)\lambda}(A|B)$ and $\tilde{T}_{\mu,k+1,\pm \lambda}(A|B)$. 

\begin{proposition}
If $\lambda > 0, \mu \in \mathbb{R}$ and $k = 0,1,2,\ldots$, then 
\begin{description}
\item[(1)]  $\displaystyle{S_{\mu-(k+1)\lambda}(A|B) \leq \tilde{T}_{\mu,k+1,-\lambda}(A|B) \leq S_{\mu-k\lambda}(A|B)}$. 
\item[(2)]  $\displaystyle{S_{\mu+k\lambda}(A|B) \leq \tilde{T}_{\mu,k+1,\lambda}(A|B) \leq S_{\mu+(k+1)\lambda}(A|B)}$.
\end{description}
\label{prop:proposition3}
\end{proposition}

\begin{flushleft}
{\bf Proof.}  If $\lambda > 0, \mu \in \mathbb{R}$ and $k = 0,1,2,\ldots$, 
then it is easy to give the following inequalities for any $t >0$ : 
\end{flushleft}
$$
t^{\mu-(k+1)\lambda}\log t \leq \frac{t^{\mu-(k+1)\lambda}-t^{\mu-k\lambda}}{-\lambda} \leq t^{\mu-k\lambda} \log t, 
$$
$$
t^{\mu+k\lambda}\log t \leq \frac{t^{\mu+(k+1)\lambda}-t^{\mu+k\lambda}}{\lambda} \leq t^{\mu+(k+1)\lambda}\log t.
$$
Then replace $t$ by $A^{-1/2}BA^{-1/2}$ and multiply $A^{1/2}$ on both sides so we get the desired results.  
\ \hfill q.e.d. 

\vspace{0.5cm}

By putting $k = 0$ or $1$,  we get the following. 

\begin{corollary}
For $A > 0, B > 0, \mu \in \mathbb{R}$ and $\lambda > 0$,
\begin{eqnarray*}
&  & S_{\mu-2\lambda}(A|B) \leq \tilde{T}_{\mu,2,-\lambda}(A|B) \leq S_{\mu-\lambda}(A|B) \\
& \leq & \tilde{T}_{\mu,1,-\lambda}(A|B) \leq S_{\mu}(A|B) \leq \tilde{T}_{\mu,1,\lambda}(A|B) \\
& \leq & S_{\mu+\lambda}(A|B) \leq \tilde{T}_{\mu,2,\lambda}(A|B) \leq S_{\mu+2\lambda}(A|B).
\end{eqnarray*}
\label{cor:corollary2}
\end{corollary}

In particular by putting $\mu = 0, \lambda = 1$, we get the following. 

\begin{corollary}
For $A > 0, B > 0$, 
\begin{eqnarray*}
&   & S_{-2}(A|B) \leq \tilde{T}_{0,2,-1}(A|B) \leq S_{-1}(A|B) \\
& \leq & \tilde{T}_{0,1,-1}(A|B) \leq S_{0}(A|B) \leq \tilde{T}_{0,1,1}(A|B) \\
& \leq & S_{1}(A|B) \leq \tilde{T}_{0,2,1}(A|B) \leq S_{2}(A|B).
\end{eqnarray*}
We rewrite the following: 
\begin{eqnarray*}
&   & S_{-2}(A|B) \leq AB^{-1}A-AB^{-1}AB^{-1}A \leq S_{-1}(A|B) \\
& \leq & A-AB^{-1}A \leq S(A|B) \leq B-A \\
& \leq & S_{1}(A|B) \leq BA^{-1}B-B \leq S_{2}(A|B).
\end{eqnarray*}
\label{cor:corollary3}
\end{corollary}

Similarly we state the relationship among 
$\sum_{j=1}^n S_{\mu \pm k\lambda}(A_j|B_j), \\
\sum_{j=1}^n S_{\mu \pm (k+1)\lambda}(A_j|B_j)$ and 
$\sum_{j=1}^n \tilde{T}_{\mu,k+1,\pm \lambda}(A_j|B_j)$, where $A_j > 0, B_j > 0$ \\
satisfying $\sum_{j=1}^n A_j = \sum_{j=1}^n B_j = I$. If 
$\lambda > 0, \mu \in \mathbb{R}$ and $k = 0,1,2,\ldots$, then we have the 
following: 
$$
\sum_{j=1}^n S_{\mu-(k+1)\lambda}(A_j|B_j) \leq \sum_{j=1}^n \tilde{T}_{\mu,k+1,-\lambda}(A_j|B_j) \leq \sum_{j=1}^n S_{\mu-k\lambda}(A_j|B_j), 
$$
$$
\sum_{j=1}^n S_{\mu+k\lambda}(A_j|B_j) \leq \sum_{j=1}^n \tilde{T}_{\mu,k+1,\lambda}(A_j|B_j) \leq \sum_{j=1}^n S_{\mu+(k+1)\lambda}(A_j|B_j). 
$$

By putting $k = 0$ or $1$, we get the following. 

\begin{corollary}
For $A > 0, B > 0, \mu \in \mathbb{R}$ and $\lambda > 0$, 
\begin{eqnarray*}
&  & \sum_{j=1}^n S_{\mu-2\lambda}(A_j|B_j) \leq \sum_{j=1}^n \tilde{T}_{\mu,2,-\lambda}(A_j|B_j) \leq \sum_{j=1}^n S_{\mu-\lambda}(A_j|B_j) \\
& \leq & \sum_{j=1}^n \tilde{T}_{\mu,1,-\lambda}(A_j|B_j) \leq \sum_{j=1}^n S_{\mu}(A_j|B_j) \leq \sum_{j=1}^n \tilde{T}_{\mu,1,\lambda}(A_j|B_j) \\
& \leq & \sum_{j=1}^n S_{\mu+\lambda}(A_j|B_j) \leq \sum_{j=1}^n \tilde{T}_{\mu,2,\lambda}(A_j|B_j) \leq \sum_{j=1}^n S_{\mu+2\lambda}(A_j|B_j).
\end{eqnarray*}
\label{cor:corollary4}
\end{corollary}

In particular by putting $\mu = 0, \lambda = 1$, we get the following result which is somewhat different type from Corollary in \cite{Fu:par}.

\begin{corollary}
For $A_j > 0, B_j > 0$ satisfying $\sum_{j=1}^n A_j = \sum_{j=1}^n B_j = I$,
\begin{eqnarray*}
&   & \sum_{j=1}^n S_{-2}(A_j|B_j) \leq \sum_{j=1}^n A_jB_j^{-1}A_j-\sum_{j=1}^n A_jB_j^{-1}A_jB_j^{-1}A_j \\
& \leq & \sum_{j=1}^n S_{-1}(A_j|B_j) \leq I-\sum_{j=1}^n A_jB_j^{-1}A_j \leq \sum_{j=1}^n S(A_j|B_j) \leq 0 \\
& \leq & \sum_{j=1}^n S_{1}(A_j|B_j) \leq \sum_{j=1}^n B_jA_j^{-1}B_j-I \leq \sum_{j=1}^n S_{2}(A_j|B_j).
\end{eqnarray*}
\label{cor:corollary5}
\end{corollary}

%
\section*{Acknowledgement}
The authours thank Professor T. Furuta for reading our manuscript and for valuable comments.

\vspace{1cm}

\end{document}